\documentclass[11pt]{article}
\usepackage{amsmath}
\usepackage{amsfonts}
\usepackage{cite}

\numberwithin{equation}{section} \topmargin=0cm
\date{}

\begin{document}

\title {\bf On certain type of difference polynomials of meromorphic functions }

\author{ Ranran Zhang$^{1}$\ , \ Zhibo Huang$^{2,}$\footnote{Corresponding
author}\\
\footnotesize
{$^{1}$Department of Mathematics, Guangdong University of Education, Guangzhou 510303,  China}\\
\footnotesize
{$^{2}$School of Mathematical Sciences, South China Normal University, Guangzhou 510631,  China}\\
\footnotesize {E-mail: $^{1}$zhangranran@gdei.edu.cn\ , \ $^{2}$huangzhibo@scnu.edu.cn}}

\date{}
\maketitle

{\small \noindent{\bf Abstract}\quad
In this paper, we investigate zeros of difference polynomials of the form
$f(z)^nH(z, f)-s(z)$, where $f(z)$ is a meromorphic function, $H(z, f)$ is a difference polynomial of $f(z)$ and $s(z)$ is
a small function.
We first obtain some inequalities for the relationship of the zero counting function of $f(z)^nH(z, f)-s(z)$ and the
characteristic function and pole counting function of $f(z)$.
Based on these inequalities, we establish some difference analogues of a classical result of Hayman for meromorphic functions.
Some special cases are also investigated.
These results improve previous findings.
\\

\\

 }

\section{ Introduction and results }

Let $f(z)$ be a meromorphic function in the complex plane
$\mathbb{C}$. We assume that the reader is familiar with the basic
notions of Nevanlinna's theory (see \cite{hay2}). We use $\sigma(f)$ to
denote the order of growth of $f(z)$, $\sigma_2(f)$ to denote the hyper order of $f(z)$, and $\delta(\infty, f)$ to denote
the Nevanlinna deficiency of $f(z)$.
Moreover, we denote by $S(r, f)$ any real function of growth $o(T(r, f))$ as
$r\rightarrow\infty$ outside of a possible exceptional set of finite
logarithmic measure. A meromorphic funtion $\alpha(z)$ is said to be a small function of $f(z)$, if $T(r, \alpha)=S(r, f)$.

Many authors have been interested in the value distribution of differential polynomials of meromorphic functions
and obtained fruitful results. In particular, Hayman proved the following results.
\vskip%
0.1in

\noindent{\bf Theorem A} (\cite{hay1})\quad \emph{ If $f(z)$ is a transcendental entire function and $n\geq 2$,
then $f'(z)f(z)^n$ assumes all  finite values except possibly zero infinitely often.}
\vskip%
0.1in
\noindent{\bf Theorem B} (\cite{hay1})\quad \emph{ If $f(z)$ is a transcendental meromorphic function and $n\geq 3$,
then $f'(z)f(z)^n$ assumes all finite values except possibly zero infinitely often.}
\vskip%
0.1in

The difference analogues of Nevanlinna value distribution theory have been established in \cite{chi, hal2, hal3, hal, lai2}. Using
these theories, many authors considered the value distribution of difference polynomials.
The results they got  are mostly about entire functions.
In particular, the following result can be viewed as a difference analogue of Theorem A.

\vskip%
0.1in
\noindent{\bf Theorem C} (\cite{lai3, liu, zheng})\quad \emph{Let $f(z)$ be a transcendental
entire function of finite order, and $c$ be a non-zero
complex constant. Then for $n\geq 2$, $f(z)^nf(z+c)$
assumes every non-zero value $a\in\mathbb{C}$ infinitely often.}
\vskip%
0.1in

For meromorphic functions, it is easy to see that a direct difference analogue of Theorem B cannot hold.
Indeed, take $f(z)=\tan z$. Then
\begin{equation*}
f(z)^3f(z+\frac{\pi}{2})=-\tan^2z
\end{equation*}
never takes the value 1. A natural question is:
What can be said about the conclusion of Theorem B
if $f'(z)$ of Theorem B is replaced by $f(z+\eta)$?
For this question, the following results are obtained in \cite{lix, liu2}.
\vskip%
0.1in

\noindent{\bf Theorem D} (\cite{lix})\quad \emph{Let $f(z)$ be a transcendental meromorphic function such
that its order $\sigma(f)<\infty$, let $\eta$ be a non-zero complex number, and
let $n \geq 1$ be an integer. Suppose that $P(z) \not\equiv 0$ is a polynomial. Then}
\begin{align*}
&\overline{N}\left(r, \frac{1}{f(z)^nf(z+\eta)-P(z)}\right)\geq nT(r, f(z))+m(r, f(z))\\
&- 2\overline{N}(r, f(z))-2\overline{N}\left(r, \frac{1}{f(z)}\right)-N\left(r, \frac{1}{f(z)}\right)+o\left(\frac{T(r, f(z))}{r^{1-\varepsilon}}\right)+O(1),
\end{align*}
\emph{as $r\not\in E$ and $r\rightarrow\infty$,
where $E$ denotes a set of finite logarithmic measure.}
\vskip%
0.1in

\noindent{\bf Theorem E} (\cite{lix, liu2})\quad \emph{Let $f(z)$ be a transcendental meromorphic function such
that its order $\sigma(f)<\infty$, let $\eta$ be a non-zero complex number, and
let $n \geq 6$ be an integer.  Suppose that $P(z) \not\equiv 0$ is a polynomial. Then
$f(z)^nf(z+\eta)-P(z)$ has infinitely many zeros.}
\vskip%
0.1in

We pose three questions related to Theorems D and E.
\vskip%
0.1in

1. What happens if $f(z+\eta)$ is generalized to difference polynomials?

2. Is it possible to reduce the
condition ``$n \geq 6$'' in Theorem E?

3. Applying Theorem D, we cannot get Theorem C. So Theorem D is not a direct improvement of Theorem C to the case of meromorphic functions.
Is it possible to obtain  such a direct improvement?
\vskip%
0.1in

In this paper, we consider these questions and obtain some results using different methods than \cite{lix, liu2}.
Among our results, Theorem 1.1 and Corollary 1.1 answer questions 1 and 3, and
Corollaries 1.2, 1.3 and 1.4 offer partial results concerning question 2.
\vskip%
0.1in

To formulate our results, we introduce some notations.
The difference
polynomial $H(z, f)$ of a meromorphic function $f(z)$ is defined by
$$
H(z, f)=\sum_{\lambda\in
J}a_\lambda(z)\prod_{j=1}^{\tau_\lambda}f(z+\delta_{\lambda,
j})^{\mu_{\lambda, j}},\eqno(1.1)
$$
where $J$ is an index set, $\delta_{\lambda, j}$ are complex
constants, $\mu_{\lambda, j}$ are non-negative integers, and the
coefficients $a_\lambda(z)(\not\equiv0)$ are small meromorphic
functions of $f(z)$. The degree of the monomial $a_\lambda(z)\prod_{j=1}^{\tau_\lambda}f(z+\delta_{\lambda,
j})^{\mu_{\lambda, j}}$ is defined by
$$
d_{\lambda}=\sum_{j=1}^{\tau_\lambda}\mu_{\lambda, j}.\eqno(1.2)
$$
The degree of $H(z, f)$ is defined by
$$
d_H=\deg_fH(z, f)=\max_{\lambda\in J}d_{\lambda}.\eqno(1.3)
$$
Let the different $\delta_{\lambda, j}$ in $H(z, f)$ be $\delta_1, \cdots,
\delta_m$, and let
$$
\chi=\left\{
\begin{array}{lll}
1,\ \mbox{if}\ \delta_s=0 \ \mbox{for some}\ s\in\{1,\cdots,m\},  \\
0,\ \mbox{if}\ \delta_t\neq0\ \mbox{for all}\ t=1, \cdots, m.
\end{array}\right.\eqno(1.4)
$$
\vskip%
0.1in

\noindent{\bf Theorem 1.1}\quad  \emph{Let $f(z)$ be a transcendental meromorphic function satisfying $\sigma_2(f)<1$, let $H(z, f)(\not\equiv0)$
be a difference polynomial in
$f(z)$ of the form (1.1) with $m\geq1$ different $\delta_{\lambda, j}$, let $d_H$ and $\chi$ be defined by (1.3) and (1.4) respectively,
and let $n>md_H$  be an integer.  If $s(z) \not\equiv 0$ is a small meromorphic function of $f(z)$, then}
\begin{align*}
&2\overline{N}\left(r, \frac{1}{f(z)^nH(z, f)-s(z)}\right)\geq(n-1)T(r, f(z))\\
&- (m-\chi)d_H N(r, f(z))-(2m+1-2\chi)\overline{N}(r, f(z))+S(r, f).
\end{align*}
\vskip%
0.1in

For a difference monomial
$$
F(z,f)=f(z+c_1)^{i_1}f(z+c_2)^{i_2}\cdots f(z+c_m)^{i_m},\eqno(1.5)
$$
where $m\geq1$ is an integer, $i_1, i_2, \cdots, i_m$ are positive integers, and $c_1, c_2, \cdots, c_m$ are different non-zero
complex constants, we obtain the following corollary.
\vskip%
0.1in

\noindent{\bf Corollary 1.1}\quad  \emph{Let $f(z)$ be a transcendental meromorphic function satisfying $\sigma_2(f)<1$, let $F(z, f)$
be a difference monomial in
$f(z)$ of the form (1.5) with $m\geq1$ different shifts, let $\deg_fF(z, f)=d_{F}$, and let $n>d_{F}$  be an integer.
If $s(z) \not\equiv 0$ is a small meromorphic function of $f(z)$, then}
\begin{align*}
2\overline{N}\left(r, \frac{1}{f(z)^nF(z, f)-s(z)}\right)&\geq (n-1)T(r, f(z))\\
&- d_FN(r, f(z))-(2m+1)\overline{N}(r, f(z))+S(r, f).
\end{align*}
\emph{Especially, if $F(z,f)=f(z+\eta)\ (\eta\in \mathbb{C}/\{0\})$, then}
$$
2\overline{N}\left(r, \frac{1}{f(z)^nf(z+\eta)-s(z)}\right)\geq (n-1)T(r, f(z))- N(r, f(z))-3\overline{N}(r, f(z))+S(r, f).
$$
\vskip%
0.1in

Theorem 1.1 and Corollary 1.1 generalize Theorem D to difference polynomials and are direct improvements of Theorem C to meromorphic functions.
Furthermore,
using Corollary 1.1 we can get Corollary 1.2, which is a version to reduce the
condition ``$n \geq 6$'' in Theorem E.
\vskip%
0.1in

\noindent{\bf Corollary 1.2}\quad  \emph{Let $f(z)$ be a transcendental meromorphic function
satisfying $\sigma_2(f)<1$ and $\delta(\infty, f(z))>\frac{1}{2}$, let $\eta$ be a non-zero complex number, and
let $n\geq 3$ be an integer.  If $s(z) \not\equiv 0$ is a small meromorphic function of $f(z)$, then
$f(z)^nf(z+\eta)-s(z)$ has infinitely many zeros.}
\vskip%
0.1in

For the difference monomial (1.5), if the poles, zeros and shifts of $f(z)$ satisfy some conditions,
we can obtain a better estimate.
\vskip%
0.1in

\noindent{\bf Theorem 1.2}\quad  \emph{Let $f(z)$ be a transcendental meromorphic function satisfying $\sigma_2(f)<1$, let $F(z, f)$
be a difference monomial in
$f(z)$ of the form (1.5) with $m\geq1$ different shifts, let $\deg_fF(z, f)=d_{F}$, and let $n>d_{F}$  be an integer. Suppose that all except for finitely many poles
$z_i$ and zeros $z_j$ of $f(z)$ satisfy $z_i-z_j\neq c_l$ $(l=1,\cdots,m)$.
If $s(z) \not\equiv 0$ is a small meromorphic function of $f(z)$, then}
$$
2\overline{N}\left(r, \frac{1}{f(z)^nF(z, f)-s(z)}\right)\geq (n-1)T(r, f(z))- (2m+1)\overline{N}(r, f(z))
+S(r, f).
$$
\vskip%
0.1in
From Theorem 1.2, we can easily get the following corollary, which
reduces the condition ``$n \geq 6$'' in another way.

\vskip%
0.1in
\noindent{\bf Corollary 1.3}\quad  \emph{Let $f(z)$ be a transcendental meromorphic function
satisfying $\sigma_2(f)<1$, let $\eta$ be a non-zero complex number, and
let $n\geq 5$ be an integer. Suppose that all except for finitely many poles
$z_i$ and zeros $z_j$ of $f(z)$ satisfy $z_i-z_j\neq \eta$. If $s(z) \not\equiv 0$ is a small meromorphic function of $f(z)$, then
$f(z)^nf(z+\eta)-s(z)$ has infinitely many zeros.}
\vskip%
0.1in

At last, we estimate the zeros of $f(z)^nH(z, f)-s(z)$ under the
assumption that $f(z)$ has two Borel exceptional values.
\vskip%
0.1in

\noindent{\bf Theorem 1.3}\quad  \emph{Let $f(z)$ be a finite order transcendental meromorphic function with two Borel exceptional
values $a, b\in \mathbb{C}\cup\{\infty\}$, let $H(z, f)(\not\equiv0)$
be a difference polynomial in
$f(z)$ of the form (1.1) with $m\geq1$ different $\delta_{\lambda, j}$, let $d_\lambda$ and $d_H$ be defined by (1.2) and (1.3) respectively, and
let $n$ be a positive integer.  Suppose that $s(z) \not\equiv 0$ is a small meromorphic function of $f(z)$.}

(i) \emph{If $a, b\in \mathbb{C}$, $a^n\sum\limits_{\lambda\in J}a_\lambda(z)a^{d_{\lambda}}-s(z)\not\equiv0$, $b^n\sum\limits_{\lambda\in J}a_\lambda(z)b^{d_{\lambda}}-s(z)\not\equiv0$ and $n>md_H$, then}
$$
N\left(r,\frac{1}{f(z)^nH(z, f)-s(z)}\right)\geq (n-md_H)T(r, f(z))+S(r, f).
$$

(ii) \emph{If $a\in \mathbb{C}, b=\infty$ and $a^n\sum\limits_{\lambda\in J}a_\lambda(z)a^{d_{\lambda}}-s(z)\not\equiv0$, then}
$$
N\left(r,\frac{1}{f(z)^nH(z, f)-s(z)}\right)\geq nT(r, f(z))+S(r, f).
$$
\vskip%
0.1in
From Theorem 1.3, we can easily get the following corollary, which reduces the condition ``$n \geq 6$''
to  ``$n \geq 2$'' for meromorphic functions with two Borel exceptional
values.

\vskip%
0.1in
\noindent{\bf Corollary 1.4}\quad  \emph{Let $f(z)$ be a finite order transcendental meromorphic function with two Borel exceptional
values $a, b\in \mathbb{C}\cup\{\infty\}$, let $\eta$ be a non-zero complex number, and
let $n\geq 2$ be an integer.  Suppose that $s(z) \not\equiv 0$ is a small meromorphic function of $f(z)$, and that one of the following two conditions holds:}

(i) \emph{ $a, b\in \mathbb{C}$, $a^{n+1}-s(z)\not\equiv0$ and $b^{n+1}-s(z)\not\equiv0$;}

(ii) \emph{ $a\in \mathbb{C}, b=\infty$ and $a^{n+1}-s(z)\not\equiv0$.}

\noindent\emph{Then $f(z)^nf(z+\eta)-s(z)$ has infinitely many zeros.}

\section{ Proofs of Theorem 1.1 and Corollary 1.1}
We need the following lemmas.
\vskip%
0.1in

\noindent{\bf Lemma 2.1} (\cite{hal})\quad \emph{ Let $f(z)$ be a non-constant meromorphic function
and $c\in\mathbb{C}$. If $\sigma_2(f)<1$ and $\varepsilon>0$, then
$$
m\left(r, \frac{f(z+c)}{f(z)}\right)=o\left(\frac{T(r, f(z))}{r^{1-\sigma_2(f)-\varepsilon}}\right)
$$
for all $r$ outside of a set of finite logarithmic measure.}
\vskip%
0.1in

By \cite[Lemma 1]{abl}, \cite[p. 66]{gol} and \cite[Lemma 8.3]{hal}, we immediately deduce the following lemma.
\vskip%
0.1in

\noindent{\bf Lemma 2.2}\quad \emph{ Let $f(z)$ be a non-constant meromorphic function of $\sigma_2(f)<1$,
and let $c\neq0$ be an arbitrary complex
number. Then}
$$
T\big(r, f(z+c)\big)=T(r, f(z))+S(r, f),
$$
$$
N\big(r, f(z+c)\big)=N(r, f(z))+S(r, f),
$$
$$
\overline{N}\big(r, f(z+c)\big)=\overline{N}(r, f(z))+S(r, f).
$$
\vskip%
0.1in
Applying logarithmic derivative lemma and Lemma 2.1 to Theorem 2.3 of \cite{lai2}, we get the following
lemma.
\vskip%
0.1in

\noindent{\bf Lemma 2.3}\quad \emph{ Let $f(z)$ be a transcendental
meromorphic solution of hyper order $\sigma_2(f)<1$ of a differential-difference equation of the form
\begin{equation*}
U(z, f)P(z, f)=Q(z, f),
\end{equation*}
where $U(z, f)$ is a difference polynomial in $f(z)$ with small meromorphic coefficients,
$P(z, f)$ and $Q(z, f)$ are differential-difference polynomials in
$f(z)$  such that the proximity functions of the coefficients of $P(z, f)$ and $Q(z, f)$ are of the type $S(r, f)$. Assume that $\deg_fU(z,f)=n$,
$\deg_fQ(z, f)\leq n$ and
$U(z, f)$ contains just
one term of maximal total degree in $f(z)$ and its shifts. Then}
\begin{equation*}
m\big(r, P(z, f)\big)=S(r,f).
\end{equation*}
\vskip%
0.1in

Using a similar proof as in \cite[Theorem 1.1]{zh} or \cite[Lemma 2]{zheng2}, we get the following lemma.
\vskip%
0.1in

\noindent{\bf Lemma 2.4}\quad \emph{Let $f(z)$ be a transcendental meromorphic function satisfying $\sigma_2(f)<1$, let $H(z, f)(\not\equiv0)$
be a difference polynomial in
$f(z)$ of the form (1.1) with $m\geq1$ different $\delta_{\lambda, j}$,
let $F(z, f)$ be a difference monomial in
$f(z)$ of the form (1.5), and let $\deg_fH(z, f)=d_H$ and $\deg_fF(z, f)=d_{F}$.
Then}
$$
T(r, H(z, f))\leq md_HT(r, f(z))+S(r, f),\eqno(2.1)
$$
$$
T(r, F(z, f))\leq d_{F}T(r, f(z))+S(r, f).\eqno(2.2)
$$
\vskip%
0.1in

{\bf  Proof of Theorem 1.1.}\quad  Set
$$
\psi(z)=f(z)^nH(z, f)-s(z).\eqno(2.3)
$$
First observe that
$\psi(z)\not\equiv0$. Indeed, if $\psi(z)\equiv 0$, then
$$
H(z, f)\equiv \frac{ s(z)}{f(z)^n}.\eqno(2.4)
$$
Since $n>md_H$, comparing the characteristic functions of both sides of (2.4) and using (2.1) of Lemma 2.4, we get a contradiction.
So $\psi(z)\not\equiv 0$.

Differentiating both sides of (2.3)
we obtain
$$
\psi'(z)=nf(z)^{n-1}f'(z)H(z, f)+f(z)^nH'(z, f)-s'(z).\eqno(2.5)
$$
Since $\psi(z)\not\equiv 0$, multiplying both sides of (2.3) by $\frac{\psi'(z)}{\psi(z)}$, we get
$$
\psi'(z)=\frac{\psi'(z)}{\psi(z)}f(z)^nH(z, f)-\frac{\psi'(z)}{\psi(z)}s(z).\eqno(2.6)
$$
Subtracting (2.5) from (2.6), we get
$$
f(z)^{n-1}E(z)=s'(z)-\frac{\psi'(z)}{\psi(z)}s(z),\eqno(2.7)
$$
where
$$
E(z)=nf'(z)H(z, f)-\frac{\psi'(z)}{\psi(z)}f(z)H(z, f)+f(z)H'(z, f).\eqno(2.8)
$$

We affirm that $E(z)\not\equiv0$. Otherwise, since $s(z)\not\equiv0$, it follows from (2.7) that
$$
\frac{\psi'(z)}{\psi(z)}=\frac{s'(z)}{s(z)},
$$
which gives $\psi(z)=C_1s(z),$ where $C_1$ is a non-zero constant. Substituting $\psi(z)=C_1s(z)$ into (2.3), we get
$$
H(z, f)=\frac{(C_1+1)s(z)}{f(z)^n}.\eqno(2.9)
$$
Similarly as in (2.4), by (2.9) and (2.1), we get a contradiction. So $E(z)\not\equiv0$.

By (2.1), we have $T(r, \psi(z))\leq (n+md_H)T(r, f(z))+S(r, f).$
So
$$
m\left(r, \frac{\psi'(z)}{\psi(z)}\right)=S(r, \psi)=S(r, f).\eqno(2.10)
$$
Applying Lemma 2.3 to equation (2.7), we have
$$
m(r, E(z))=S(r, f).\eqno(2.11)
$$

Next we estimate $N(r, E(z))$. By (2.8), we see that the poles of $E(z)$ come from the poles of $f(z)$, the poles of $H(z, f)$, and the poles of $\frac{\psi'(z)}{\psi(z)}$.
We denote by $N(r, |E(z)=f(z)=\infty)$
the counting function of those common poles of $E(z)$ and $f(z)$ in
$|z|<r$, where each such point is counted according to its multiplicity in $E(z)$,  denote by $N(r, |E(z)=H(z, f)=\infty, f(z)\neq \infty)$
the counting function of those common poles of $E(z)$ and $H(z, f)$ in
$|z|<r$, where each such point is not a pole of $f(z)$, and each such point is counted according to its multiplicity in $E(z)$, and
denote by $N(r, |E(z)=\frac{\psi'(z)}{\psi(z)}=\infty, f(z)\neq \infty, H(z, f)\neq \infty)$
the counting function of those common poles of $E(z)$ and $\frac{\psi'(z)}{\psi(z)}$ in
$|z|<r$, where each such point is not a pole of $f(z)$ or a pole of $H(z, f)$, and each such point is counted according to its multiplicity in $E(z)$.
Then
\begin{align}
\setcounter{equation}{11}
N(r, E(z))&=N(r, |E(z)=f(z)=\infty)\nonumber\\
&+N(r, |E(z)=H(z, f)=\infty, f(z)\neq \infty)\nonumber\\
&+N\left(r, |E(z)=\frac{\psi'(z)}{\psi(z)}=\infty, f(z)\neq \infty, H(z, f)\neq \infty\right).
\end{align}

Suppose that $z_0$ is a pole of $E(z)$ with order $k$.

If $z_0$ is a pole of $f(z)$ with order $p$,  by (2.7), $n\geq2$ and the fact that $\frac{\psi'(z)}{\psi(z)}$ has at most simple poles, we see that
$z_0$ must be a pole of $s(z)$ with order $q$ and $k+(n-1)p\leq q+1$. We then deduce from $n\geq2$ that $k\leq q$.   So
$$
N(r, |E(z)=f(z)=\infty)\leq N(r,s(z))=S(r, f).\eqno(2.13)
$$

If $z_0$ is not a pole of $f(z)$ and $z_0$ is a pole of $H(z, f)$ with order $l$, then by (2.8) and the fact that $\frac{\psi'(z)}{\psi(z)}$ has at most simple poles,
we see that $k\leq l+1$. We
denote by $N(r, |H(z,f)=\infty, f(z)\neq\infty)$ the counting function of those poles of $H(z, f)$ in
$|z|<r$, where each such point is not a pole of $f(z)$, and each such point is counted according to its multiplicity in $H(z, f)$,
and denote by $\overline{N}(r,|H(z,f)=\infty, f(z)\neq\infty)$ the counting function of those poles of $H(z,f)$ in
$|z|<r$, where each such point is not a pole of $f(z)$, and each such point is counted one time. Then
\begin{align}
\setcounter{equation}{13}
&N(r, |E(z)=H(z, f)=\infty, f(z)\neq \infty)\nonumber\\
&\leq N(r, |H(z,f)=\infty, f(z)\neq\infty)+\overline{N}(r,|H(z,f)=\infty, f(z)\neq\infty).
\end{align}
Since the
different $\delta_{\lambda, j}$ in $H(z, f)$ be $\delta_1, \cdots,
\delta_m$ and $\chi$ is defined by (1.4), we deduce from Lemma 2.2 that
\begin{align}
\setcounter{equation}{14}
N(r, |H(z,f)=\infty, f(z)\neq\infty)&\leq\sum_{j=1}^m d_H N(r, f(z+\delta_j))-\chi d_H N(r, f(z))+S(r, f) \nonumber\\
&= (m-\chi )d_HN(r, f(z))+S(r, f),
\end{align}
\begin{align}
\setcounter{equation}{15}
\overline{N}(r, |H(z,f)=\infty, f(z)\neq\infty)&\leq\sum_{j=1}^m \overline{N}\big(r, f(z+\delta_j)\big)-\chi \overline{N}(r, f(z))+S(r, f)\nonumber\\
&= (m-\chi)\overline{N}(r, f(z))+S(r, f).
\end{align}

If $z_0$ is not a pole of $f(z)$ and $z_0$ is not a pole of $H(z, f)$, then $z_0$ must be a pole of $\frac{\psi'(z)}{\psi(z)}$.
Since $\frac{\psi'(z)}{\psi(z)}$ has at most simple poles, we deduce from (2.8) that $k=1$.
The poles of $\frac{\psi'(z)}{\psi(z)}$ come from the poles of $\psi(z)$ and the zeros of $\psi(z)$.
If $z_0$ is a pole of $\psi(z)$, then by (2.3), we see that $z_0$ must be a pole of $s(z)$.
So
\begin{align}
\setcounter{equation}{16}
&N\left(r, |E(z)=\frac{\psi'(z)}{\psi(z)}=\infty, f(z)\neq \infty, H(z, f)\neq \infty\right)\nonumber\\
&\leq \overline{N}(r,s(z))+\overline{N}\left(r,\frac{1}{\psi(z)}\right)\nonumber\\
&=\overline{N}\left(r,\frac{1}{\psi(z)}\right)+S(r, f).
\end{align}
We deduce from (2.11)--(2.17) that
$$
T(r, E(z))\leq (m-\chi )d_HN(r, f(z))+(m-\chi)\overline{N}(r, f(z))+\overline{N}\left(r,\frac{1}{\psi(z)}\right)+S(r, f).\eqno(2.18)
$$

By (2.7) and (2.10), we get
\begin{align}
\setcounter{equation}{18}
(n-1)T(r, f(z))
&\leq T(r, E(z))+T\left(r,\frac{\psi'(z)}{\psi(z)}\right)+S(r, f)\nonumber\\
&=T(r, E(z))+N\left(r,\frac{\psi'(z)}{\psi(z)}\right)+S(r, f)\nonumber\\
&=T(r, E(z))+\overline{N}\left(r,\frac{1}{\psi(z)}\right)+\overline{N}(r,\psi(z))+S(r, f).
\end{align}
Since $H(z, f)$ has $m$ different $\delta_{\lambda, j}$ and $\chi$ is defined by (1.4), we deduce from (2.3) and Lemma 2.2 that
\begin{align}
\setcounter{equation}{19}
\overline{N}(r,\psi(z))&\leq\overline{N}(r, f(z))+\overline{N}(r, H(z, f))-\chi\overline{N}(r, f(z))+S(r, f)\nonumber\\
&\leq(1+m-\chi)\overline{N}(r, f(z))+S(r, f).
\end{align}
We deduce from (2.18)--(2.20) that
\begin{align*}
2\overline{N}\left(r, \frac{1}{f(z)^nH(z, f)-s(z)}\right)&=2\overline{N}\left(r,\frac{1}{\psi(z)}\right)\geq(n-1)T(r, f(z))\\
&- (m-\chi)d_H N(r, f(z))-(2m+1-2\chi)\overline{N}(r, f(z))+S(r, f).
\end{align*}

\vskip%
0.1in

{\bf  Proof of Corollary 1.1.}\quad
By (2.2) of Lemma 2.4 and using the similar method as in the proof of Theorem 1.1, we can prove Corollary 1.1 easily.

\section{ Proof of Theorem 1.2}
Set
$$
\psi(z)=f(z)^nF(z, f)-s(z).
$$
Since $n>d_{F}$, we deduce from (2.2) of Lemma 2.4 that $\psi(z)\not\equiv0$. Since $F(z, f)$ is a special case of $H(z, f)$,
we also have (2.5)--(2.12), where $H(z, f)$ is replaced by $F(z, f)$.
Next we discuss each term in (2.12).

Suppose that $z_0$ is a pole of $E(z)$ with order $k$.

If $z_0$ is a pole of $f(z)$, as in (2.13) of Theorem 1.1, we get
$$
N(r, |E(z)=f(z)=\infty)\leq N(r,s(z))=S(r, f).\eqno(3.1)
$$

If $z_0$ is not a pole of $f(z)$ and $z_0$ is a pole of $F(z, f)$, then $z_0$ must be a pole of $f(z+c_t)$ for some $t\in \{1,\cdots,m\}$. So $z_0+c_t$ is a pole of $f(z)$.
Since all except for finitely many poles
$z_i$ and zeros $z_j$ of $f(z)$ satisfy $z_i-z_j\neq c_l$ $(l=1,\cdots,m)$, we will assume that $z_0$ is not a zero of $f(z)$. So,
when estimating $N(r, |E(z)=F(z, f)=\infty, f(z)\neq \infty)$,
we may have an error term of the type $O(\log r)$.
Since $f(z_0)\neq0, \infty$, we see that $z_0$ is a pole of $f(z)^{n-1}E(z)$
with order $k$. Furthermore, $\frac{\psi'(z)}{\psi(z)}$ has at most simple poles.
By (2.7), we see that $z_0$ is a pole of $\frac{\psi'(z)}{\psi(z)}$ with order 1 and
$k=1$, or $z_0$ is a pole of $s(z)$ with order $q$ and $k\leq q+1$. So,
$z_0$ is a simple pole of $E(z)$ or $z_0$ is a pole of $s(z)$ with $k\leq q+1$.
Therefore,
$$
N(r, |E(z)=F(z, f)=\infty, f(z)\neq \infty)\leq \overline{N}(r, F(z,f))+N(r, s(z))+\overline{N}(r, s(z))+O(\log r).
$$
By Lemma 2.2, we have
$$
\overline{N}(r, F(z,f))\leq\sum_{j=1}^m \overline{N}(r, f(z+c_j))= m \overline{N}(r, f(z))+S(r, f).\eqno(3.2)
$$

If $z_0$ is not a pole of $f(z)$ and $z_0$ is not a pole of $F(z, f)$, then $z_0$ must be a pole of $\frac{\psi'(z)}{\psi(z)}$.
As in (2.17) of Theorem 1.1, we get
$$
N\left(r, |E(z)=\frac{\psi'(z)}{\psi(z)}=\infty, f(z)\neq \infty, F(z, f)\neq \infty\right)\leq \overline{N}\left(r,\frac{1}{\psi(z)}\right)+S(r, f).\eqno(3.3)
$$

By (2.11), (2.12) and (3.1)--(3.3), we get
$$
T(r, E(z))\leq m \overline{N}(r, f(z))+\overline{N}\left(r,\frac{1}{\psi(z)}\right)+S(r, f).\eqno(3.4)
$$

By (2.7) and (2.10), we get
\begin{align}
\setcounter{equation}{4}
(n-1)T(r, f(z))
\leq T(r, E(z))+\overline{N}\left(r,\frac{1}{\psi(z)}\right)+\overline{N}(r,\psi(z))+S(r, f).
\end{align}
Since $F(z, f)=f(z+c_1)^{i_1}\cdots f(z+c_m)^{i_m}$, $\psi(z)=f(z)^nF(z, f)-s(z)$ and $c_1, \cdots, c_m$ are different non-zero
complex constants, we deduce from Lemma 2.2 that
\begin{align}
\setcounter{equation}{5}
\overline{N}(r,\psi(z))&\leq\overline{N}(r, f(z))+\overline{N}(r, F(z, f))+S(r, f)\nonumber\\
&\leq(1+m)\overline{N}(r, f(z))+S(r, f).
\end{align}
We deduce from (3.4)--(3.6) that
\begin{align*}
2\overline{N}\left(r, \frac{1}{f(z)^nF(z, f)-s(z)}\right)&=2\overline{N}\left(r,\frac{1}{\psi(z)}\right)\\
&\geq (n-1)T(r, f(z))- (2m+1)\overline{N}(r, f(z))
+S(r, f).
\end{align*}

\section{ Proof of Theorem 1.3}
We need the following lemma.
\vskip%
0.1in

\noindent{\bf Lemma 4.1} (\cite{lip})\quad\emph{ Suppose that $h$ is a
non-constant meromorphic function satisfying
$$
\overline{N}(r, h)+\overline{N}(r, 1/h)=S(r, h).
$$
Let $f=a_0h^p+a_1h^{p-1}+\cdots+a_p$, and $g=b_0h^q
+b_1h^{q-1}+\cdots+b_q$ be polynomials in $h$ with coefficients
$a_0, a_1, \cdots, a_p, b_0, b_1, \cdots, b_q$ being small functions
of $h$ and $a_0b_0a_p\not\equiv0$. If $q \leq p$, then $m(r, g/f ) =
S(r,h)$.}
\vskip%
0.1in

{\bf  Proof of Theorem 1.3.}\quad
Set
$$
\psi(z)=f(z)^nH(z, f)-s(z).\eqno(4.1)
$$

First we assume that the condition (i) in Theorem 1.3 holds.
Let
$$
g(z)=\frac{f(z)-a}{f(z)-b}.
$$
Then $0, \infty$ are two Borel exceptional values of
$g(z)$. By Hadamard
factorization theorem, $g(z)$ takes the form
$$
g(z)=w(z)e^{h(z)},
$$
where $w(z)$ is a meromorphic function such that
$\sigma(w(z))<\sigma(g(z))$, and $h(z)$ is a polynomial such that $\sigma
(g(z))=\deg h(z)\geq1$. So
$$
f(z)=\frac{bw(z)e^{h(z)}-a}{w(z)e^{h(z)}-1}.\eqno(4.2)
$$

Substituting (4.2) into $f(z)^n$, we get
$$
f(z)^n=\frac{b^nw(z)^ne^{nh(z)}+\cdots+(-a)^n}{w(z)^ne^{nh(z)}+\cdots+(-1)^n}.\eqno(4.3)
$$
Denoting
$$
W_\lambda(z)=w(z+\delta_{\lambda,1})^{\mu_{\lambda, 1}}\cdots w(z+\delta_{\lambda,\tau_\lambda})^{\mu_{\lambda, \tau_\lambda}}
$$ and substituting (4.2) into $H(z, f)$, we get
\begin{align}
\setcounter{equation}{3}
&H(z, f)=\sum_{\lambda\in
J}a_\lambda(z)\prod_{j=1}^{\tau_\lambda}\frac{b^{\mu_{\lambda, j}}w(z+\delta_{\lambda,j})^{\mu_{\lambda, j}}e^{\mu_{\lambda, j}h(z+\delta_{\lambda,j})}+\cdots+(-a)^{\mu_{\lambda, j}}}
{w(z+\delta_{\lambda,j})^{\mu_{\lambda, j}}e^{\mu_{\lambda, j}h(z+\delta_{\lambda,j})}+\cdots+(-1)^{\mu_{\lambda, j}}}\nonumber\\
&=\sum_{\lambda\in
J}a_\lambda(z)\frac{b^{\mu_{\lambda, 1}+\cdots+\mu_{\lambda, \tau_\lambda}}W_\lambda(z)
e^{\mu_{\lambda, 1}h(z+\delta_{\lambda,1})+\cdots+\mu_{\lambda, \tau_\lambda}h(z+\delta_{\lambda,\tau_\lambda})}+\cdots+(-a)^{\mu_{\lambda, 1}+\cdots+\mu_{\lambda, \tau_\lambda}}}
{W_\lambda(z)
e^{\mu_{\lambda, 1}h(z+\delta_{\lambda,1})+\cdots+\mu_{\lambda, \tau_\lambda}h(z+\delta_{\lambda,\tau_\lambda})}+\cdots+(-1)^{\mu_{\lambda, 1}+\cdots+\mu_{\lambda, \tau_\lambda}}}.\nonumber
\end{align}
Denoting
$$
s_\lambda(z)=W_\lambda(z)e^{\mu_{\lambda, 1}(h(z+\delta_{\lambda,1})-h(z))}\cdots e^{\mu_{\lambda, \tau_\lambda}(h(z+\delta_{\lambda,\tau_\lambda})-h(z))},
$$
we have
\begin{align}
\setcounter{equation}{3}
&H(z, f)=\sum_{\lambda\in
J}a_\lambda(z)\frac{b^{d_{\lambda}}s_\lambda(z)e^{d_{\lambda}h(z)}+\cdots+(-a)^{d_{\lambda}}}{s_\lambda(z)e^{d_{\lambda}h(z)}+\cdots+(-1)^{d_{\lambda}}}\nonumber\\
&=\frac{(\sum\limits_{\lambda\in J}a_\lambda(z)b^{d_{\lambda}})\prod\limits_{\lambda\in
J}s_\lambda(z)e^{\sum\limits_{\lambda\in J}d_{\lambda}h(z)}+\cdots
+(\sum\limits_{\lambda\in J}a_\lambda(z)a^{d_{\lambda}})(-1)^{\sum\limits_{\lambda\in J}d_{\lambda}}}
{\prod\limits_{\lambda\in
J}s_\lambda(z)e^{\sum\limits_{\lambda\in J}d_{\lambda}h(z)}+\cdots+(-1)^{\sum\limits_{\lambda\in J}d_{\lambda}}}.
\end{align}
By (4.1), (4.3), (4.4) and denoting $S(z)=\prod\limits_{\lambda\in
J}s_\lambda(z)$, $D=\sum\limits_{\lambda\in J}d_{\lambda}$, we get
\begin{align}
\setcounter{equation}{4}
&\psi(z)=\frac{b^nw(z)^n(\sum\limits_{\lambda\in J}a_\lambda(z)b^{d_{\lambda}})S(z)e^{(n+D)h(z)}+\cdots
+(\sum\limits_{\lambda\in J}a_\lambda(z)a^{d_{\lambda}})(-a)^n(-1)^{D}}
{w(z)^nS(z)e^{(n+D)h(z)}+\cdots+(-1)^{n+D}}-s(z)\nonumber\\
&=\frac{(b^n\sum\limits_{\lambda\in J}a_\lambda(z)b^{d_{\lambda}}-s(z))w(z)^nS(z)e^{(n+D)h(z)}+\cdots
+(-1)^{n+D}(a^n\sum\limits_{\lambda\in J}a_\lambda(z)a^{d_{\lambda}}-s(z))}
{w(z)^nS(z)e^{(n+D)h(z)}+\cdots+(-1)^{n+D}}.
\end{align}
We see that
$\psi(z)$ is a rational
function in $e^{h(z)}$ and
the coefficients in (4.5) are all small
functions of $e^{h(z)}$. Since $a^n\sum\limits_{\lambda\in J}a_\lambda(z)a^{d_{\lambda}}-s(z)\not\equiv0$ and $b^n\sum\limits_{\lambda\in J}a_\lambda(z)b^{d_{\lambda}}-s(z)\not\equiv0$,
by Lemma 4.1, we get
$$
m\left(r, \frac{1}{\psi(z)}\right)=S(r, e^{h(z)})=S(r, f).
$$
Moreover, by (4.1) and Lemma 2.4, we have
$$
nT(r, f(z))=T\left(r, \frac{\psi(z)+s(z)}{H(z, f)}\right)\leq T(r, \psi(z))+md_HT(r, f(z))+S(r, f).
$$
So
$$
N\left(r,\frac{1}{f(z)^nH(z, f)-s(z)}\right)=N\left(r,\frac{1}{\psi(z)}\right)\geq (n-md_H)T(r, f(z))+S(r, f).
$$

Now we assume that the condition (ii) in Theorem 1.3 holds.
Then $f(z)$ takes the form
$$
f(z)=w(z)e^{h(z)}+a,\eqno(4.6)
$$
where $w(z)$ is a meromorphic function such that
$\sigma(w(z))<\sigma(f(z))$, and $h(z)$ is a polynomial such that $\sigma
(f(z))=\deg h(z)\geq1$.
Substituting (4.6) into $f(z)^n$, we get
$$
f(z)^n=w(z)^ne^{nh(z)}+\cdots+a^n.\eqno(4.7)
$$
Using the notations $W_\lambda(z)$ and $s_\lambda(z)$ as above and substituting (4.6) into $H(z, f)$, we get
\begin{align*}
H(z, f)&=\sum_{\lambda\in J}a_\lambda(z)\prod_{j=1}^{\tau_\lambda}(w(z+\delta_{\lambda,j})^{\mu_{\lambda, j}}e^{\mu_{\lambda, j}h(z+\delta_{\lambda,j})}+\cdots+a^{\mu_{\lambda, j}})\\
&=\sum_{\lambda\in J}a_\lambda(z)(W_\lambda(z)
e^{\mu_{\lambda, 1}h(z+\delta_{\lambda,1})+\cdots+\mu_{\lambda, \tau_\lambda}h(z+\delta_{\lambda,\tau_\lambda})}+\cdots+a^{\mu_{\lambda, 1}+\cdots+\mu_{\lambda, \tau_\lambda}})\\
&=\sum_{\lambda\in J}a_\lambda(z)(s_\lambda(z)e^{d_{\lambda}h(z)}+\cdots+a^{d_{\lambda}}).
\end{align*}
Since $H(z, f)\not\equiv0$ and  $d_H=\max\limits_{\lambda\in J}d_\lambda$,  we see that $H(z, f)$ takes the form
$$
H(z, f)=\sum_{\lambda\in J}a_\lambda(z)a^{d_{\lambda}}\not\equiv0,\eqno(4.8)
$$
or
$$
H(z, f)=l_q(z)e^{qh(z)}+\cdots+l_1(z)e^{h(z)}+\sum_{\lambda\in J}a_\lambda(z)a^{d_{\lambda}}\ (1\leq q\leq d_H), \eqno(4.9)
$$
where  $l_j(z)(j=1, \cdots, q)$ are all small functions of $e^{h(z)}$ and $l_q(z)\not\equiv0$.

If (4.8) holds, by (4.1), (4.7) and Lemma 4.1, we get
$$
N\left(r,\frac{1}{\psi(z)}\right)=nT(r, f(z))+S(r, f).
$$
If (4.9) holds, by (4.1), (4.7) and Lemma 4.1, we get
$$
N\left(r,\frac{1}{\psi(z)}\right)=(n+q)T(r, f(z))+S(r, f) \ (1\leq q\leq d_H).
$$
Therefore,
$$
N\left(r,\frac{1}{f(z)^nH(z, f)-s(z)}\right)=N\left(r,\frac{1}{\psi(z)}\right)\geq nT(r, f(z))+S(r, f).
$$
\vskip.2cm\par\quad
 \\
 {\bf Acknowledgements}\\
This work is supported by Guangdong National Natural Science Foundation of China (Nos. 2016A030313745, 2014A030313422) and Training Plan Fund of Outstanding
Young Teachers of Higher Learning Institutions of Guangdong Province of China (No. Yq20145084602).

\end{document}